\documentclass[10pt]{article}
\usepackage{amsfonts,amsmath,amssymb,amsthm}
\usepackage{enumerate,footnote,indentfirst,latexsym,times,tikz,tkz-berge}
\usepackage[top=1.9in,bottom=1.9in,left=2.0in,right=2.0in]{geometry}
\usepackage{caption,subcaption}
\usepackage[hypertexnames=true]{hyperref}
\usepackage{color}

\newcommand{\ben}{\begin{enumerate}}
\newcommand{\een}{\end{enumerate}}
\newcommand{\ble}{\begin{lem}}
\newcommand{\ele}{\end{lem}}
\newcommand{\bth}{\begin{thm}}
\renewcommand{\eth}{\end{thm}}
\newcommand{\bpr}{\begin{prop}}
\newcommand{\epr}{\end{prop}}
\newcommand{\bco}{\begin{cor}}
\newcommand{\eco}{\end{cor}}
\newcommand{\bcon}{\begin{conj}}
\newcommand{\econ}{\end{conj}}
\newcommand{\bde}{\begin{defn}}
\newcommand{\ede}{\end{defn}}
\newcommand{\bex}{\begin{exa}}
\newcommand{\eex}{\end{exa}}
\newcommand{\barr}{\begin{array}}
\newcommand{\earr}{\end{array}}
\newcommand{\btab}{\begin{tabular}}
\newcommand{\etab}{\end{tabular}}
\newcommand{\beq}{\begin{equation}}
\newcommand{\eeq}{\end{equation}}
\newcommand{\bea}{\begin{eqnarray*}}
\newcommand{\eea}{\end{eqnarray*}}
\newcommand{\bal}{\begin{align*}}
\newcommand{\bce}{\begin{center}}
\newcommand{\ece}{\end{center}}
\newcommand{\bpi}{\begin{picture}}
\newcommand{\epi}{\end{picture}}
\newcommand{\bpp}{\begin{picture}}
\newcommand{\epp}{\end{picture}}
\newcommand{\bfi}{\begin{figure} \begin{center}}
\newcommand{\efi}{\end{center} \end{figure}}
\newcommand{\bprf}{\begin{proof}}
\newcommand{\eprf}{\end{proof}\medskip}
\newcommand{\capt}{\caption}
\newcommand{\bsl}{\begin{slide}{}}
\newcommand{\esl}{\end{slide}}
\newcommand{\bfr}{\begin{frame}}
\newcommand{\efr}{\end{frame}}

\newcommand{\ol}{\overline}
\newcommand{\hor}{\mbox{--}}

\newcommand{\hso}[1]{\hspace{-1pt}}
\newcommand{\vs}[1]{\vspace{#1}}

\def\<{\langle}
\def\>{\rangle}

\newcommand{\ree}[1]{(\ref{#1})}

\newcommand{\Kb}{\ol{K}}

\hypersetup{colorlinks,citecolor=black,filecolor=black,linkcolor=black,urlcolor=black}
\newtheorem{thm}{Theorem}[section]
\newtheorem{prop}[thm]{Proposition}
\newtheorem{cor}[thm]{Corollary}
\newtheorem{lem}[thm]{Lemma}
\newtheorem{conj}[thm]{Conjecture}
\newtheorem{exa}[thm]{Example}

\makeindex

\begin{document}

\title{On Constructing Regular\\Distance-Preserving Graphs}
\author{Dennis Ross$^1$, Bruce Sagan$^2$, Ronald Nussbaum$^1$,\\and Abdol-Hossein Esfahanian$^1$\\
Michigan State University\\
East Lansing, MI 48824, U.S.A.\\}
\date{}
\maketitle

\footnotetext[1]{Department of Computer Science and Engineering,
$\{$\href{mailto:rossdenn@cse.msu.edu}{rossdenn}, \href{mailto:ronald@cse.msu.edu}{ronald}, \href{mailto:esfahanian@cse.msu.edu}{esfahanian}$\}$@cse.msu.edu}
\footnotetext[2]{Department of Mathematics,
sagan@math.msu.edu}

\begin{abstract}
Let $G$ be a simple, connected graph on $n$ vertices. Let $d_G(u,v)$ denote the distance between vertices $u$ and $v$ in $G$.
A subgraph $H$ of $G$ is \textit{isometric} if $d_H(u,v)=d_G(u,v)$ for every $u,v\in V(H)$.
We say that $G$ is a \textit{distance-preserving graph} if $G$ contains at least one isometric subgraph of order $k$ for every $k$, $1\le k\le n$.
In this paper we construct regular distance-preserving graphs of all possible orders and degrees of regularity.
By modifying the Havel-Hakimi algorithm, we are able to construct distance preserving graphs for certain other degree sequences as well.  We include a discussion of some related conjectures which we have computationally verified for small values of $n$.
\end{abstract}

\section{Introduction}
\label{sec:introduction}

In this paper we will only consider connected graphs unless noted otherwise.
A subgraph $H$ of a graph $G=(V,E)$ is \textit{isometric} if, for every pair of vertices $u$ and $v$ of $H$, we have $d_H(u,v)=d_G(u,v)$, where $d$ denotes distance.
If every connected induced subgraph of a graph is isometric, then the graph is said to be \textit{distance-hereditary}.
Originally mentioned by Sachs \cite{sachs1970berge} while working with perfect graphs, this class of graphs was first named and characterized by Howorka \cite{howorka1977characterization}.
Distance-hereditary graphs have been studied extensively in the literature (see, for example, the articles \cite{bandelt1986distance,hammer1990completely}).
They may be recognized in linear time \cite{damiand2001simple}.

We call a connected graph \textit{distance-preserving} (dp) if it contains an isometric subgraph of every possible order.
In previous papers we began to characterize dp graphs \cite{nussbaum2012preliminary} and explore potential applications \cite{nussbaum2010clustering,nussbaum2013clustering}.
In the next section we will construct $r$-regular dp graphs on $n$ vertices for all possible values of $n$ and $r$ where such graphs exist.
It was conjectured by Nussbaum and Esfahanian \cite{nussbaum2012preliminary} that almost all graphs are dp.
The results in Table \ref{tab:percentage} lead us to make a similar conjecture for regular graphs.

\begin{table}
  \center
  \begin{tabular}{|l|c|c|c|}
    \hline
    & \# connected & \# connected & \\
    $n$ & regular graphs & regular dp graphs & \% dp graphs\\
    \hline
    \hline
	5 & 2 & 1 & 50.000\\
    \hline
	6 & 5 & 4 & 80.000\\
    \hline
	7 & 4 & 3 & 75.000\\
    \hline
	8 & 17 & 14 & 82.353\\
    \hline
	9 & 22 & 20 & 90.909\\
    \hline
	10 & 167 & 153 & 91.617\\
    \hline
	11 &  539 & 484 & 89.796\\
    \hline
	12 &  18979 & 18405 & 96.976\\
    \hline
	13 &  389436 & 384319 & 98.686\\
    \hline
  \end{tabular}
  \caption{Percentage of regular graphs which are distance preserving}
  \label{tab:percentage}
\end{table}

\bcon
\label{con:almost_all_regular}
Almost all connected regular graphs are dp.
\econ

Rather than looking only at regular graphs, one could try to construct a dp graph for every possible graphical degree sequence.
The Havel-Hakimi algorithm \cite{hakimi1962realizability,hakimi1963realizability} can be used to generate a graph having a given degree sequence, although it is easy to come up with examples where the resulting graph is not dp.
However, we show in Section~\ref{sec:arbitrary_ds} that a slight modification of the algorithm can be used to produce dp graphs for certain degree sequences.

\section{Constructing regular  dp graphs}
\label{sec:regular_dp}

We will construct, for each possible $n$ and $r$, an $r$-regular dp graph with $n$ nodes.
In particular this means that, in order to have vertices of degree $r$, we must have $n\ge r+1$.  And if $r$ is odd, then $n$ must be even.
Finally, we must have $r\ge3$ since if $r=2$ then a connected graph of this regularity must be a cycle which is not dp for $n\ge5$.  We will call the remaining pairs {\em admissible}.

One class of $r$-regular $n$-vertex (but not necessarily dp) graphs which we will need are the {\em circulant graphs} $C_{n,r}$.  These can be constructed by using a  vertex set $V=\{1,\dots,n\}$ with all edges $ij$ such that $1\le |i-j| \le r/2$, together with all edges of the form $i(i+n/2)$ if $r$ is odd.  All arithmetic is being done modulo $n$. 

We will also be using two operations to construct larger graphs from smaller ones.  The {\em join} of $G$ and $H$, $G+H$, is obtained from their disjoint union by adding all edges of the form $uv$ where $u\in V(G)$ and $v\in V(H)$.  
In this case we will refer to the {\em underlying bipartite graph} of $G+H$ which consists of all the added edges.

The other operation will depend on an edge $ux$ being selected from $E(G)$ and another edge $vy$ being selected from $E(H)$.  In that case we define the {\em direct sum} of $G$ and $H$, $G\oplus H$, to be the graph obtained from their disjoint union by removing $ux$ and $vy$ and adding new edges $uv$ and $xy$.  
See the first graph  in Figure~\ref{oplus} for an example.
Our notation does not specify which edges are to be used.  But in the graphs we will use, this choice can usually be made arbitrarily because of symmetry.  If this is not the case, an explicit choice will be made.

\bfi
\begin{tikzpicture}
\draw(-1.5,1) node {$K_5\oplus K_5=$};
\fill(0,0) circle(.1);
\fill(0,1) circle(.1);
\fill(1,2) circle(.1);
\fill(2,0) circle(.1);
\fill(2,1) circle(.1);
\fill(3,0) circle(.1);
\fill(3,1) circle(.1);
\fill(4,2) circle(.1);
\fill(5,0) circle(.1);
\fill(5,1) circle(.1);
\draw(2,-.5) node{$x$};
\draw(2,1.5) node{$u$};
\draw(3,-.5) node{$y$};
\draw(3,1.5) node{$v$};
\draw (0,0)--(2,0)--(0,1)--(2,1)--(1,2)--(0,1)--(0,0) (2,0)--(1,2)--(0,0)--(2,1)--(3,1)--(5,0)--(4,2)--(3,0)--(2,0)
(5,0)--(3,0)--(5,1)--(3,1)--(4,2)--(5,1)--(5,0);
\end{tikzpicture}

\vs{10pt}

\begin{tikzpicture}
\draw(-1,1) node {$K_5^{\oplus 3}=$};
\fill(0,0) circle(.1);
\fill(0,1) circle(.1);
\fill(1,2) circle(.1);
\fill(2,0) circle(.1);
\fill(2,1) circle(.1);
\fill(3,0) circle(.1);
\fill(3,1) circle(.1);
\fill(4,2) circle(.1);
\fill(5,0) circle(.1);
\fill(5,1) circle(.1);
\fill(6,0) circle(.1);
\fill(6,1) circle(.1);
\fill(7,2) circle(.1);
\fill(8,0) circle(.1);
\fill(8,1) circle(.1);
\draw(2,-.5) node{$r$};
\draw(2,1.5) node{$a$};
\draw(3,-.5) node{$s$};
\draw(3,1.5) node{$b$};
\draw(5,-.5) node{$x$};
\draw(5,1.5) node{$u$};
\draw(6,-.5) node{$y$};
\draw(6,1.5) node{$v$};
\draw (0,0)--(2,0)--(0,1)--(2,1)--(1,2)--(0,1)--(0,0) (2,0)--(1,2)--(0,0)--(2,1)--(3,1)--(5,0)--(4,2)--(3,0)--(2,0)
(5,0)--(3,0)--(5,1)--(3,1)--(4,2)--(5,1) (3,1)--(6,1)--(8,0)--(7,2)--(6,0)--(5,0) (8,0)--(6,0)--(8,1)--(6,1)--(7,2)--(8,1)--(8,0)
;
\end{tikzpicture}
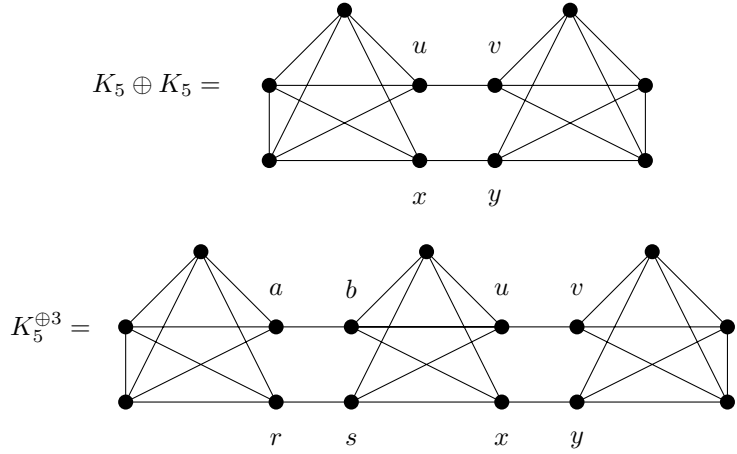
\capt{The graphs $K_5\oplus K_5$ and $K_5^{\oplus 3}$ 
\label{oplus}}
\efi

We want to extend the direct sum to more than two graphs.  So, for example,  $G\oplus H\oplus K$ will mean that we have chosen an edge $ar$ is $G$, two edges $bs$ and $ux$ in $H$, and an edge $vy$ in $K$.  These edges will be removed and replaced by $ab$, $rs$, $uv$, and $xy$.  Our only restriction is that the edges $bs$ and $ux$ are independent in $H$, that is, do not share a vertex.  
The extension to more than three graphs is done in the obvious way and we denote by $G^{\oplus k}$ the $k$-fold direct sum of $G$ with itself.
Figure~\ref{oplus} shows $K_5\oplus K_5\oplus K_5=K_5^{\oplus 3}$.

Finally, we will need the following lemma.  The ideas behind the proof are due to Zahedi \cite{zahedi2014distance}.  We use $N(v)$ to denote the neighborhood of $v$, that is, all vertices adjacent to $v$.
\ble
\label{zah}
Let $G=(V,E)$ be a graph and $v\in V$.  Suppose that  any two nonadjacent vertices $x,y\in N(v)$ have a common neighbor other than $v$.  Then $G-v$ is an isometric subgraph of $G$.
\ele
\bprf
It suffices to show that removing $v$ from $G$ does not destroy all geodesics between any two vertices $r,s\in V-v$.  So let $P$ be any $r\hor s$ geodesic in $G$.  If $P$ does not contain $v$ then we are done.  Otherwise, the vertices on $P$ just before and after $v$ are some $x,y\in N(v)$.  We can not have $xy\in E$ since then there would be a shorter $r\hor s$ path.  Thus, by assumption, there is some common neighbor $w\neq v$ of $x,y$.  It follows that $P-v+w$ is an $r\hor s$ geodesic in $G-v$ which completes the proof.
\eprf

\bth
\label{theorem:regular_construction}
For each admissible pair $(n,r)$, there exists a dp graph with $n$ vertices which is regular of degree $r$.
\eth
\bprf
Fix $r$.  We will denote the graph on $n$ vertices which we construct by $G_n$.  We will have a special construction for $r=3$, so assume for the moment that $r\ge 4$.  

For $r+1\le n\le 2r$, we let $G_n=C_{r,2r-n}+\Kb_{n-r}$ where $\Kb_s$ is the complement of the complete graph on $s$ vertices.  
The graph $G_7$ for $r=4$ is depicted on the left in Figure~\ref{Gn}.
It is easy to check that $G_n$ has $n$ vertices and is regular of degree $r$.  Also, Lemma~\ref{zah} and the completeness of the  underlying bipartite graph show that removing any vertex from either of the parts $C_{r,2r-n}$ or $\Kb_{n-r}$ gives an isometric subgraph as long as the part from which it  is  removed has at least two vertices. Continuing to delete vertices in this fashion shows that $G_n$ is dp.

\bfi
\begin{tikzpicture}
\draw (-.5,.5) node {$G_7=$};
\fill(0,0) circle(.1);
\fill(1,0) circle(.1);
\fill(1,1) circle(.1);
\fill(2,0) circle(.1);
\fill(2,1) circle(.1);
\fill(3,0) circle(.1);
\fill(3,1) circle(.1);
\foreach \a in {0,1,2,3} \foreach \b in {1,2,3}
	\draw (\a,0)--(\b,1);
\draw (0,0)--(1,0) (2,0)--(3,0);
\end{tikzpicture}
\qquad
\qquad
\begin{tikzpicture}
\draw (-1,.85) node{$x$};
\draw (-1.75,.5) node{$G_9=$};
\fill(-1,.5) circle(.1);
\fill(0,0) circle(.1);
\fill(0,1) circle(.1);
\fill(1,0) circle(.1);
\fill(1,1) circle(.1);
\fill(2,0) circle(.1);
\fill(2,1) circle(.1);
\fill(3,0) circle(.1);
\fill(3,1) circle(.1);
\foreach \a in {0,1,2,3} \foreach \b in {2,3}
	\draw (\a,0)--(\b,1);
\foreach \a in {2,3} \foreach \b in {0,1,2,3}
	\draw (\a,0)--(\b,1);
\draw (-1,.5)--(0,0)--(1,1)--(-1,.5)--(0,1)--(1,0)--(-1,.5);
\end{tikzpicture}
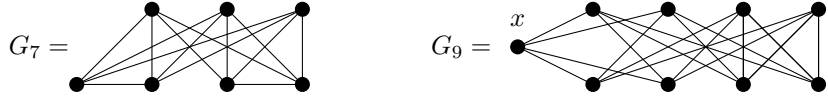
\capt{The graphs $G_7$ on the left and $G_9$ on the right for $r=4$
\label{Gn}}
\efi

To construct $G_{2r+1}$ first note that, since $n=2r+1$ is odd, we only need such a construction for $r$ even.
Take the complete bipartite graph $K_{r,r}$ and remove $r/2$ independent edges.  Now add a new vertex adjacent to precisely those vertices where an edge was removed.  This vertex will be called the {\em external vertex}, $x$.  
See Figure~\ref{Gn} for an example when $r=4$.
To construct the isometric subgraphs, first remove one by one all but two of the vertices adjacent to $x$, where those two vertices are in different parts of the underlying bipartite graph $K_{r,r}$.  These subgraphs will be isometric by Lemma~\ref{zah} and the existence of the two remaining vertices adjacent to $x$.  Now continue to remove two more vertices  not adjacent to $x$, one from each part of $K_{r,r}$.  Since $r\ge4$, there will still be a vertex of $K_{r,r}$ not adjacent to $x$ remaining in both parts and so Lemma~\ref{zah} can be applied.  For the set of $r+1$ vertices to remove, we take $x\cup N(x)$.  This leaves the graph $K_{r/2,r/2}$ which is isometric.  From here, one continues as in the $n\le 2r$ case.

To form $G_n$ for $n\ge 2r+2$, we take advantage of the fact that $K_{r+1}$ is regular of degree $r$.  So divide $n$ by $r+1$ to obtain $n=q(r+1)+t$ where $q\ge2$ and $0\le t\le r$.  So 
$$
n=(q-1)(r+1)+(t+r+1)=p(r+1)+s
$$
where $p=q-1\ge1$ and  $s=t+r+1$ so that $r+1\le s\le 2r+1$.  Now let $G_n=G_s\oplus K_{r+1}^{\oplus p}$ where the edge removed from $G_s$ is one of the edges of the underlying bipartite graph.  If $s=2r+1$, then we choose this edge to be one whose endpoints have degree $r$ in the underlying bipartite graph.  To set up notation, let $H_0$ be the subgraph of $G_n$ obtained from $G_s$ and $H_1,\dots H_p$ be the subgraphs obtained from the copies of $K_{r+1}$ listed so that $H_i$ is adjacent to $H_{i-1}$ and $H_{i+1}$ for $1\le i<r$.  We also let $ux$, $vy$, and $wz$ be  edges removed to form $H_0$, $H_1$, and $H_p$ respectively, where the edges added between $H_0$ and $H_1$ are $uv$ and $xy$.  (There may be a second edge removed to form $H_1$, but we will not need a notation for it.)

To construct the isometric subgraphs we remove vertices from $H_0$ in the same way as we did for $G_s$, taking care never to remove $u$ or $x$, until only $4$ vertices remain: $u$, $x$, and another vertex from each part of the underlying complete subgraph of $G_s$.  The resulting subgraphs will be isometric in $G_n$ because all geodesics from a vertex of $H_0$ to another $H_i$ go through $u$ or $x$ which are never removed.  Next we continue to remove two more vertices from $H_p$, neither of them being $w,z$.  Since $r\ge4$, these subgraphs will again be isometric since there will always exist a third vertex vertex of $H_p$ present in the subgraph.  So we have isometric subgraphs of order up to $s-2$.  To remove $s-1$ vertices, we delete every vertex of $H_0$ except $u$.  This is again isometric since no geodesic from $u$ to $G_n-H_0$ goes through any other vertex of $H_0$.  Next we remove $u$, which clearly leaves an isometric subgraph.  We can now continue to remove vertices in $H_1,H_2,\dots,H_p$ in that order using the same ideas as for removing the vertices of $H_0$ and always removing the vertices of $H_i$ which were adjacent to $H_{i-1}$ first.  This will produce an isometric subgraph of every order.

\bfi
\begin{tikzpicture}
\draw(0,.5) node {$G_{12}=$};
\fill(1,0) circle(.1);
\fill(1,1) circle(.1);
\fill(2,0) circle(.1);
\fill(2,1) circle(.1);
\fill(3,0) circle(.1);
\fill(3,1) circle(.1);
\fill(4,0) circle(.1);
\fill(4,1) circle(.1);
\fill(5,0) circle(.1);
\fill(5,1) circle(.1);
\fill(6,0) circle(.1);
\fill(6,1) circle(.1);
\foreach \x in {1,2,3,4,5} \foreach \y in {0,1}
	\draw (\x,\y)--(\x+1,\y);
\foreach \x in {1,3,4,6}
	\draw (\x,0)--(\x,1);
\draw (1,0)--(2,1) (1,1)--(2,0) (5,0)--(6,1) (6,0)--(5,1);
\foreach \x in {1,2,3,4,5,6}
	\draw (\x,-.5) node{$v_{\x}$};
\foreach \x in {1,2,3,4,5,6}
	\draw (\x,1.5) node{$u_{\x}$};
\end{tikzpicture}
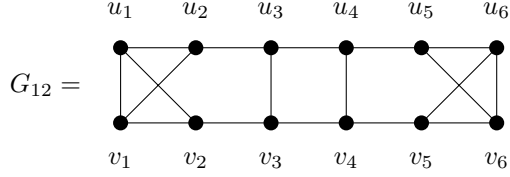
\capt{The graph $G_{12}$ when $r=3$
\label{r=3}}
\efi

There remains to do the case $r=3$.  Since $r$ is odd, $n$ must be even so let $k=n/2$.  Let $G$ have vertices $u_1,\dots,u_k$ and $v_1,\dots,v_k$ where the two sets of vertices form paths in their given orders.  Also add  all edges of the form $u_j v_j$ except for $u_2,v_2$ and $u_{k-1} v_{k-1}$.  Finally add the edges  $u_1 v_2, u_2v_1, u_{k-1} v_k, u_k v_{k-1}$. 
The graph $G_{12}$ will be found in Figure~\ref{r=3}.
 We now just list the sets of vertices to be removed.  Verifying that the resulting subgraphs are isometric is routine.  We start with
$$
\{u_1\},\{u_1,u_k\},\{u_1,v_1,u_2\},\{u_1,v_1,u_2,u_k\},\{u_1,v_1,u_2,v_2,u_k\},
$$
and
$$
\{u_1,v_1,u_2,v_{k-1},u_k,v_k\}.
$$
We continue by deleting the $6$ vertices in the last set and adding to them the following vertices taken sequentially in the given order
$$
u_3,u_4,\dots,u_{k-1},v_2,v_3,\dots,v_{k-2}.
$$
This completes the proof.
\eprf

\section{Arbitrary graphical degree sequences}
\label{sec:arbitrary_ds}

We would like to construct a dp graph for any graphical integer sequence for which such a graph exists.
In Theorem~\ref{th:hakimi_havel} we prove that a modified version of the Havel-Hakimi algorithm generates a dp graph when no reordering of the vertex degrees is done, even when such a reordering would be called for by the original algorithm.

Let us recall the usual Havel-Hakimi algorithm.  The input is a weakly decreasing integral sequence $(d_1,\dots,d_n)$, and the output is a graph with this degree sequence, if one exists, in which case the sequence is called \textit{graphical}.
The main loop of the algorithm is as follows, where the vertices of the graph we are trying to construct are $v_1,\dots,v_n$.
\ben
\item[(a)] Add edges from $v_1$ to $v_2,\dots,v_{d_1+1}$.
\item[(b)] Remove $d_1$ from the sequence and subtract one from $d_2,\dots,d_{d_1+1}$.
\item[(c)] Remove any resulting zeros in the new sequence and rearrange the rest to be weakly decreasing, updating indices as necessary.
\een
One iterates this loop until one of two possible outcomes result.  If the algorithm reaches the empty sequence  then it has constructed a graph $G$ with $\deg v_i=d_i$ for all $i$ and this is  called a \textit{successful termination}.  Otherwise, the algorithm halts because it becomes impossible to perform step (a) in which case the sequence is not graphical, called an \textit{unsuccessful termination}.

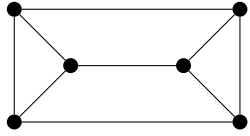
\begin{figure}
  \begin{center}
    \begin{tikzpicture}
      \SetVertexNoLabel
      \tikzstyle{VertexStyle} = [shape=circle, fill=black, minimum size=2pt,inner sep=2pt]
      \tikzstyle{EdgeStyle} = [thin]
      \Vertex[x=0, y=0]{0}
      \Vertex[x=3, y=0]{1}
      \Vertex[x=.75, y=.75]{2}
      \Vertex[x=2.25, y=.75]{3}
      \Vertex[x=0, y=1.5]{4}
      \Vertex[x=3, y=1.5]{5}
      \Edge(0)(1)
      \Edge(0)(2)
      \Edge(0)(4)
      \Edge(1)(3)
      \Edge(1)(5)
      \Edge(2)(3)
      \Edge(2)(4)
      \Edge(3)(5)
      \Edge(4)(5)
    \end{tikzpicture}
  \end{center}
  \caption{A graphical realization of $(3,3,3,3,3,3)$}
  \label{fig:hh_non_dp}
\end{figure}

In the modified version of this algorithm which we will consider, one does not rearrange the degree sequence in step (c).  This will mean that the new algorithm may have an unsuccessful termination on a sequence which is graphical.
For example, the sequence $(3, 3, 3, 3, 3, 3)$ is clearly graphical, as seen in Figure \ref{fig:hh_non_dp}.
However, three iterations of our modified version of the Havel-Hakimi algorithm leave us with $(3, 3)$, and the algorithm terminates unsuccessfully.
On the positive side, we will show that when the procedure does terminate successfully with a connected graph $G$, then $G$ must be dp.  In fact, it will turn out that the graphs
\beq
\label{Gi}
G_i = \text{subgraph of $G$ induced by $v_1,\dots,v_i$}
\eeq
will be the desired isometric subgraphs.

\bth
\label{th:hakimi_havel}
Let $S=(d_1,\ldots,d_n)$ be a weakly decreasing sequence of  integers on which the modified Havel-Hakimi algorithm terminates successfully with a connected graph $G$.  Then $G$ is dp with isometric subgraphs $G_i$ as in~\ree{Gi}.
\eth

\bprf
We will induct on $n$, where the result is obvious for $n=1$.  So assume the result for sequences with $n-1$ elements and let $S=(d_1,\dots,d_n)$ be such that the modified algorithm terminates successfully with a connected graph $G$.
To apply induction, we must first show that the algorithm terminates successfully on the degree sequence for the graph $G'=G_{n-1}$.  Applying the procedure to $G$ and $G'$ will be exactly the same until one comes to some vertex $v_i$ which is adjacent to $v_n$ in $G$.  But in the sequence $S'$ for $G'$ we will have $\deg_{G'} v_i=\deg_G v_i -1$.  So when processing $S'$, $v_i$ will be attached to exactly the same vertices as for $S$, with the exception of $v_n$.  It follows that since $S$ was brought to a successful conclusion, so must $S'$.

We must also show that $G'$ is connected.  Since $G$ is connected and $G'=G-v_n$, it suffices to show that $v_n$ is not a cut vertex.  We will actually prove the stronger statement that any two neighbors of $v_n$ are adjacent.  So take $v_i,v_j\in N(v_n)$ where we can assume, without loss of generality, that $i<j$.  Since we do not reorder vertices and $v_i$ is adjacent to $v_n$, it must also be adjacent to all vertices which have subscript greater than $i$  and whose modified degree is still positive when $v_i$ is processed by the main loop.  Since $v_j$ is such a vertex, we have proved our claim.

Now we can apply induction so that $G_1,\dots,G_{n-1}=G'$ are isometric subgraphs of $G'$.  Thus we will be done if we can show that $G'$ is isometric in $G$.  For this, it suffices to show that no geodesic of $G$ goes through $v_n$.  But, as proved in the previous paragraph, any two neighbors of $v_n$ are adjacent.  It follows that any path $P$ through $v_n$ can be made shorter by replacing the edges into and out of $v_n$ by the single edge between its neighbors on $P$.  Hence $P$ is not a geodesic and the proof is complete.
\eprf

\begin{table}
  \center
  \begin{tabular}{|l|c|c|c|}
    \hline
    & \# graphical & &\\
    $n$ & degree sequences & \# successes & \% successes\\
    \hline
    \hline
	5 & 20 & 12 & 60.000\\
    \hline
	6 & 71 & 32 & 45.070\\
    \hline
	7 & 240 & 86 & 35.833\\
    \hline
	8 & 871 & 243 & 27.899\\
    \hline
	9 & 3148 & 703 & 22.332\\
    \hline
	10 & 11655 & 2094 & 17.967\\
    \hline
	11 & 43332 & 6369 & 14.698\\
    \hline
	12 & 162769 & 19770 & 12.146\\
    \hline
  \end{tabular}
  \caption{Success rate of the modified Havel-Hakimi algorithm}
  \label{tab:hh_mod}
\end{table}

In Table \ref{tab:hh_mod} the results of the modified Havel-Hakimi algorithm are provided for $5\leq n\leq 12$.
Note that any time the algorithm terminates properly is counted as a success, whether the resulting graph is connected or not.
So, unfortunately, it seems as if our previous result only applies to a vanishingly small percentage of graphs.
It would be very interesting to find a method which would produce a dp graph for a larger class of degree sequences, especially since this might result in some progress on the conjectures mentioned in the introduction.\\

\textit{Acknowledgment .
We thank Brendan McKay, creator of the nauty software package \cite{mckay2014practical} which was used in our computations.}\\


\begin{thebibliography}{10}

\bibitem{bandelt1986distance}
H.J. Bandelt and H.M. Mulder.
\newblock {Distance-Hereditary Graphs}.
\newblock {\em Journal of Combinatorial Theory, Series B}, 41(2):182--208,
  1986.

\bibitem{damiand2001simple}
G.~Damiand, M.~Habib, and C.~Paul.
\newblock {A simple paradigm for graph recognition: application to cographs and
  distance hereditary graphs}.
\newblock {\em Theoretical Computer Science}, 263(1-2):99--111, 2001.

\bibitem{hakimi1962realizability}
SL~Hakimi.
\newblock On realizability of a set of integers as degrees of the vertices of a
  linear graph. i.
\newblock {\em Journal of the Society for Industrial \& Applied Mathematics},
  10(3):496--506, 1962.

\bibitem{hakimi1963realizability}
SL~Hakimi.
\newblock On realizability of a set of integers as degrees of the vertices of a
  linear graph ii. uniqueness.
\newblock {\em Journal of the Society for Industrial \& Applied Mathematics},
  11(1):135--147, 1963.

\bibitem{hammer1990completely}
P.L. Hammer and F.~Maffray.
\newblock {Completely Separable Graphs*}.
\newblock {\em Discrete Applied Mathematics}, 27(1-2):85--99, 1990.

\bibitem{howorka1977characterization}
E.~Howorka.
\newblock {A characterization of distance-hereditary graphs}.
\newblock {\em The Quarterly Journal of Mathematics}, 28(4):417, 1977.

\bibitem{mckay2014practical}
B.D. McKay and A.~Piperno.
\newblock Practical graph isomorphism, ii.
\newblock {\em Journal of Symbolic Computation}, 60:94--112, 2014.

\bibitem{nussbaum2012preliminary}
R.~Nussbaum and A-H. Esfahanian.
\newblock Preliminary results on distance-preserving graphs.
\newblock {\em Combinatorial Structures and their Applications, Gordon and
  Beach, New York}, 37:384, 1970.

\bibitem{nussbaum2010clustering}
R.~Nussbaum, A-H. Esfahanian, and P-N. Tan.
\newblock {Clustering Social Networks Using Distance-Preserving Subgraphs}.
\newblock In {\em 2010 International Conference on Advances in Social Networks
  Analysis and Mining (ASONAM)}, pages 380--385. IEEE, 2010.

\bibitem{nussbaum2013clustering}
R.~Nussbaum, A-H. Esfahanian, and P-N. Tan.
\newblock {Clustering Social Networks Using Distance-preserving Subgraphs}.
\newblock In T.~{\"O}zyer, editor, {\em The Influence of Technology on Social
  Network Analysis and Mining}, pages 331--349. Springer, 2013.

\bibitem{sachs1970berge}
H.~Sachs.
\newblock On the berge conjecture concerning perfect graphs.
\newblock {\em Combinatorial Structures and their Applications, Gordon and
  Beach, New York}, 37:384, 1970.

\bibitem{zahedi2014distance}
E.~Zahedi.
\newblock Distance preserving graphs.
\newblock Preprint.

\end{thebibliography}
\end{document}